\documentclass{article}
\usepackage{psfrag}
\usepackage{amssymb,bbm}
\usepackage{amsmath}
\usepackage{amsfonts}
\usepackage{pstricks}
\usepackage{graphicx}
\usepackage{yfonts}
\usepackage{setspace}
\newcounter{theorem}

\newtheorem{proposition}{Proposition}
\newtheorem{remark}{Remark}
\newtheorem{definition}{Definition}
\newtheorem{corollary}{Corollary}
\newenvironment{proof}[1][Proof]{\textbf{#1.} }{\rule{0.5em}{0.5em}}

\begin{document}

\title{The Critical Point of a Sigmoidal Curve:  the  Generalized Logistic Equation Example}
\markboth{The Critical Point of a Sigmoidal Curve}{The Critical Point of a Sigmoidal Curve}

\author{Ayse Humeyra Bilge and Yunus Ozdemir}

\maketitle

\begin{abstract}
\baselineskip=12pt \noindent
Let $y(t)$ be a smooth sigmoidal curve,  $y^{(n)}(t)$ be its  $n$th derivative, $\{t_{m,i}\}$ and  $\{t_{a,i}\}$, $i=1,2,\dots$  be the set of points where respectively the derivatives of odd and  even order reach their extreme values. The ``critical point of the sigmoidal curve" is defined to be the common limit of the  sequences $\{t_{m,i}\}$ and $\{t_{a,i}\}$, provided that the limit exists.  We prove that if $f(t)=\frac{dy}{dt}$ is an even function such that the magnitude of the analytic representation $|f_A(t)|=|f(t)+if_h(t)|$, where $f_h(t)$ is the Hilbert transform of $f(t)$, is monotone on $(0,\infty)$, then the point $t=0$ is the critical point in the sense above.  For the  general case, where $f(t)=\frac{dy}{dx}$ is not even, we prove that if $|f_A(t)|$ monotone on $(0,\infty)$ and if the phase of its Fourier transform $F(\omega)$ has a limit as $\omega\to \pm\infty$, then $t=0$ is still the critical point but as opposed to the previous case, the maximum of $f(t)$  is located away from $t=0$. We compute the Fourier transform of the generalized logistic growth functions and illustrate the notions above on  these examples.
\end{abstract}



\section{Introduction}
A sigmoidal function $y(t)$ is a monotone increasing function with horizontal asymptotes as $t\to \pm \infty$.
In the course of a study of the gelation phenomena \cite{BP2012}, we have observed an interesting  property of sigmoidal functions representing solutions of epidemic models:  The points $\{t_n\}$ where the  derivatives $y^{(n)}(t)$ reach their global extreme value seemed  to form a convergent sequence. This limit point was interpreted as a
junction point representing the critical point of a phase transition.
 In \cite{BP2013}, we have used a combination of symbolic and numeric computations  to compute and display the derivatives of various sigmoidal functions up to orders $100-200$.   For sigmoidal curves with an even derivative, we have {\it observed} that the global extrema of the odd derivatives occur at $t=0$, while the global extrema of the even derivatives tend   towards $t=0$.  As examples of sigmoidal curves without symmetry, we have studied the generalized logistic growth and the Susceptible-Infected-Removed epidemic model. For such curves, we have again observed that the global extrema of the even and odd derivatives tend  to a point in between the second and third derivatives, but the rates of   ``convergence" of the  even and odd derivatives are quite different.  We have also  studied the dependency of the location of the critical point to the system parameters and we have seen that it  moves towards the zero of the third derivative of the sigmoidal function as the asymmetry of the curve increases, making a sharper turn before its inflection point \cite{BP2013}.
 Despite strong evidence of the existence of the limit, we were unable to give a proof of the existence of a limit point, even in the case of simplest examples.

In the present work we use Fourier and Hilbert transform techniques to prove the existence of the critical point of a sigmoidal curve under certain assumptions.
The motivation for the choice of our assumptions is the following observation: Let $f(t)$ be the first derivative of a sigmoidal curve and $F(\omega)$ be its Fourier transform.
If $f(t)$ and all its derivatives are localized functions in time, then as the order of differentiation increases, $f^{(n)}(t)$ looks like a wave packet in the time domain. In the frequency domain,  the  Fourier transform of $f^{(n)}$, $F^{(n)}(\omega)=(i \omega)^n F(\omega)$ appears as a shift of $F(\omega)$ to increasingly higher frequencies $\pm \omega_n$.  If $f(t)$ is symmetrical, then $F(\omega)$ is real and the envelope of the inverse transform is always located at $t=0$. This proves the existence of the critical point in the symmetrical case.
In the general case, $F(\omega)$ has  a nonzero phase. We recall that a linear phase factor in $F(\omega)$ corresponds to a time shift in $f(t)$.
 Since the main frequency contribution to  $f^{(n)}$ is near $\pm\omega_n$, we look at the the phase of $F(\omega)$
  at $\pm\omega_n$.  If this phase $\phi_n$ is bounded or grows at most linearly with $\omega$ as $n\to\infty$, then $\phi_n/\omega_n$, hence the time shift in the derivatives will have a limit hence there is a critical point. By an appropriate  choice of the origin in the time domain we can locate the critical point at $t=0$ as in the symmetrical case.

 The definitions and theorems necessary for subsequent derivations are presented in Section 2.  In Section 3, we first prove that if $y(t)$ is an odd sigmoidal function then  $t=0$ is the critical point, provided that  certain conditions on the Hilbert transform of its derivatives  are satisfied. Then we show that if the sigmoidal function has no symmetry, then the existence of the critical point is proven with an additional hypothesis on the phase of its Fourier transform.  The computation of the Fourier transform of the generalized logistic family is presented in Appendix.   In Section 4, we use this result to illustrate how the assumptions of Section 3 are satisfied for the generalized logistic family.

\section{Preliminaries}

\subsection{Sigmoidal curves}

Let $y(t) $ be a monotone increasing function with horizontal asymptotes $y_1$ and $y_2$ as $t\to \pm \infty$ and with $\displaystyle \lim_{t\to \pm \infty}y^{(n)}(t)=0$ for all $n\ge 1$.
The first derivative $y'(t)$ has a horizontal tangent at least at one point say $t_m$,  the second derivative $y''(t)$ vanishes at $t_m$, hence  $y'''(t)$ has horizontal tangents at, at least two points $t_a$ and $t_b$ at the right and at the left of $t_m$. It follows that the third derivative $y'''(t)$ has at least $2$ zeros.   By similar arguments we can conclude  that the $n$th derivative $y^{(n)}(t)$ has at least $n-1$ zeros.  A sigmoidal curve $y(t)$ where $y^{(n)}(t)$ has exactly $n-1$ zeros is, in a sense,  a curve which makes no ``unnecessary oscillations".
 We will call such a curve ``a basic sigmoidal curve" \cite{BP2013}.

\begin{definition}\label{def1-basicsigmoidalcurve}
 A basic sigmoidal curve is a curve $y(t)$ such that
 \[\lim_{t\to -\infty} y(t)=y_1,\quad \lim_{t\to \infty}y(t)=y_2, \quad y_1\ne y_2, \quad \lim_{t\to \pm\infty} y^{(n)}(t)=0,\quad  {\rm for \ all\ } n\]
 and $y^{(n)}(t)$ has exactly $n-1$ zeros.
\end{definition}

We propose the following definition for the ``critical point".

\begin{definition}\label{def2-criticalpoint}
Let $y(t)$ be a sigmoidal curve and assume that the set of points where the even derivatives and the odd derivatives reach
their absolute extremum converge to the same point. The common limit of these derivatives, if it exists, is called the  critical
point of the sigmoidal curve.
\end{definition}

\begin{remark} \label{remark1}
The curve $y(t)=e^{-t^2}(1+t^2)$ is not basic sigmoidal. For the first few $n$,  the number of zeros of the derivatives $y^{(n)}$ are greater than $n-1$ but this sigmoidal function has a critical point.  There are strong indications that  if there exists $n_0 \in \mathbb{N}$ such that the derivatives $y^{(n)}$ form an Intrinsic Mode Function family  for $n\ge n_0$, then the sigmoidal function has a critical point. But such generalizations are not addressed here. (IMF is a function that has exactly one zero in between any two extrema \cite{IMA1} and note that derivatives of a basic sigmoidal curve form an IMF  family.)
\end{remark}

\subsection{The Fourier and Hilbert transforms}

We start by giving the definition of the Fourier and Hilbert transforms  and state their  properties.

\vskip 0.2cm
\noindent
{\bf The Fourier transform: }
The Fourier transform of a function $f(t)$, ${\cal F}(f)=F(\omega)$ is defined as
$$F(\omega)=\frac{1}{\sqrt{2\pi}}\ \int_{-\infty}^{\infty}
\ f(t) \, e^{-i\omega t}\ dt,$$
provided that the integral exists in the sense of Cauchy principal value \cite{PapoulisFourier}.
If $f(t)$ is in $L^1$, then its Fourier transform exists.
Since a sigmoidal function is finite as $t\to\infty$, its first derivative is in $L^1$.
We can recover $f(t)$ from the inverse transform by
$$f(t)=\frac{1}{\sqrt{2\pi}}\ \int_{-\infty}^{\infty}
\ F(\omega)\, e^{i\omega t}\ d\omega.$$

\vskip 0.2cm
\noindent
{\bf The Hilbert transform and the analytic representation:}
For our purposes, the simplest description of the Hilbert transform is given by its relation to the Fourier transform \cite{Duoan}.
Given $f(t)$ and its Fourier transform $F(\omega)$, we define the function $F_A(\omega)$ by

$$F_A(\omega)=\left\{
                 \begin{array}{ccl}
                   2F(\omega) & , & \quad \omega>0 \\
                   0          & , & \quad \omega<0
                 \end{array}
               \right. .
$$
The inverse transform is a complex function that is called the ``analytic representation" $f_A(t)$ of $f(t)$. The imaginary part of $f_A(t)$ is the Hilbert transform $f_h(t)$ of $f(t)$.
$$f_A(t)=f(t)+if_h(t)=\sqrt{\frac{2}{\pi}} \int_0^\infty F(\omega) \ e^{i\omega t}\ d\omega.$$

\vskip 0.2cm
\noindent
{\bf Symmetry properties:}
The property, ${\cal F}\left(\overline{f(t)}\right) =\overline{F(-\omega)}$ implies that if $f(t)$ is real then
$F(-\omega)=\overline{F(\omega)}$.
Thus if $f(t)$ is real,  $F(\omega)$ will be real provided that
$F(-\omega)=F(\omega)$.
The scaling  property
$f( at)\to  \frac{1}{|a|} F\left(\frac{\omega}{a}\right) $
implies that
$f(-t) \to  F(-\omega)$,
hence if $f(t)$ is real and even, then $F(\omega)$ is real and even.
Similarly, if $f(t)$ is real and odd, then $F(\omega)$ is pure imaginary and odd.
We note that if $f(t)$ is even (odd), its Hilbert transform is odd (even).

\vskip 0.2cm
\noindent
{\bf Differentiation:}
The effect of differentiation in the time domain is multiplication by $i\omega$ in the frequency domain.
Thus
$$f^{(n)}(t)\to (i\omega)^n F(\omega).$$

\vskip 0.2cm
\noindent
{\bf Convolution and modulation:}
There is a correspondence between products and convolutions in the time and frequency domains;
 multiplication in the time domain leads to convolution in the frequency domain, i.e,
$$f(t) \, g(t) \to \frac{1}{\sqrt{2\pi}} \, F(\omega) \, G(\omega).$$
The ``modulation" of a low frequency signal in the time domain is the multiplication of this signal  by a sinusoidal function of fixed (usually high) angular frequency $\omega_0$. In the frequency domain, the Fourier transform of the low frequency function is convolved with the Fourier transform of the sinusoid. The Fourier transform of a  pure sinusoid is not defined in the usual sense, but it is represented as the Dirac $\delta$ functions occurring at $\pm \omega_0$ and convolution carries the spectrum of the low frequency signal to the frequencies  $\pm \omega_0$.
Since the Fourier transform of a complex exponential is a $\delta$-function, we have the  correspondence below:
$$f(t)e^{i\omega_0t}\to F(\omega-\omega_0).$$

\vskip 0.2cm
\noindent
{\bf Time shift:}
As an analogue of multiplication with a complex exponential in the time domain, the  multiplication of a function in the  frequency domain by a linear phase factor leads to a shift in the time domain:
$$e^{-i\alpha \omega}F(\omega)\to f(t-\alpha).$$
\vskip 0.2cm

\noindent
{\bf Separating the even and odd components  of a function $f(t)$:}
Let $f(t)$ be the first derivative of an odd sigmoidal function. Note that even if $f(t)$  is  even, a translation of its graph to $t_o$ will be no longer even,
hence the Fourier transform will be no longer real. Nevertheless, we can recognize this situation because the phase will be linear in $\omega$.
If $\tilde{f}(t)$ is a perturbation of the original function, its Fourier transform will be complex.  If the phase of its Fourier transform remains bounded as $\omega\to \pm \infty$, then we can
distinguish between phase factors resulting from a time shift and the ones that are intrinsically related to the asymmetry of the  shape of the curve  in the time domain.  This leads us to an invariant  definition of the even/odd  components  of a function.

\begin{definition}\label{def3-evencomponent}
 Let $f(t)$ be a  function such that the Fourier transform $F(\omega)$ exists and the phase of $F(\omega)$ remains bounded as $\omega\to\pm \infty$.  The even component of $f(t)$ is the inverse Fourier transform of the real part of $F(\omega)$.
\end{definition}

\vskip 0.2cm
\noindent
{\bf The envelope:}
The magnitude of the analytic representation $f_A(t)$ gives the envelope of $f(t)$ in a sense as remarked below.

\begin{remark}\label{remark-01-alterne-roots}
Writing $f_A(t)=A(t)e^{i\varphi(t)},$
we can express $f(t)$ and $f_h(t)$ as
$$f(t)=\left|f_A(t)\right| \ \cos\left( \varphi(t)\right),\quad f_h(t)=A(t)\sin\left( \varphi(t)\right).$$
If $t_1$ and $t_2$ are two consecutive zeros of $f(t)$, then we should have
$\varphi(t_1)=\frac{\pi}{2}+k\pi$ and $\varphi(t_2)=\frac{\pi}{2}+(k+1)\pi$.  Thus, provided that $\varphi(t)$ is continuous,
there will be a time
$t_3$, $t_1<t_3<t_2$, such that $\varphi(t_3)=(k+1)\pi$, hence, $f_h(t_3)=0.$  Thus the zeros of $f(t)$ and $f_h(t)$ alternate.
It follows that $f(t_3)=\pm \left|f_A(t_3)\right|$ for some $t_1<t_3<t_2$.
\end{remark}

\section{Critical Points of Sigmoidal Curves}

Let $y(t)$ be a sigmoidal curve and let $f(t)=y^{(1)}(t)$.  Starting from $f(t)$, we can compute its Fourier transform $F(\omega)$, the Fourier transform of its analytic representation $F_A(\omega)$, the analytic representation $f_A(t)$ and the magnitude of the analytical representation $|f_A(t)|$:
$$f(t)\longrightarrow F(\omega)\longrightarrow F_A(\omega)
\longrightarrow f_A(t)\longrightarrow |f_A(t)|.$$
The same sequence of operations is applied to the derivatives of $f(t)$ as
$$f^{(n)}(t)\rightarrow F^{(n)}(\omega)=(i\omega)^n F(\omega) \rightarrow F_A(\omega)= (i\omega)^n F_A(\omega)
\rightarrow f^{(n)}_A(t)\rightarrow |f^{(n)}_A(t)|.$$
Loosely speaking, if $F(\omega)$ is localized and the effect of multiplication by $\omega^n$ is a shift of $F(\omega)$ to  some frequency $\omega_n$, then $f^{(n)}(t)$ will be approximately a modulation of its envelope, given by $|f_A(t)|$ multiplied by a sinusoid of angular frequency $\omega_n$.
If the phase of $F(\omega)$ at $\omega_n$ is $\phi_n$, then $\pm f^{(2n)}(t)$ will be a shift of the modulated waveform to $\alpha_{2n}=\phi_{2n}/\omega_{2n}$. For odd $n$, the same scheme applies to the Hilbert transform of $f^{(n)}(t)$.  If the phase of $F(\omega)$ remains bounded as $\omega\to \pm \infty$, then the time shift $\alpha_n$ will approach zero.  For even $f(t)$, all even derivatives and the Hilbert transform of the odd derivatives will be wave packets centered at $t=0$. For the asymmetrical case they will be wave packets that move towards $t=0$ as the degree of differentiation tends to infinity.  Thus the critical point will always be $t=0$, provided that the original curve is appropriately shifted in time, so that the phase of its Fourier transform is bounded as $\omega\to\pm \infty$.  In this section we will prove these assertions first for odd sigmoidal functions, then for the general case.

\subsection{Odd sigmoidal curves}
We start with the symmetrical case.
Let $y(t)$ be a sigmoidal function whose first derivative $y^{(1)}(t)$ is even. Then, $y(t)$ up to an additive constant will be odd and successive derivatives will have similar properties.  The Fourier transform of the odd sigmoidal function exists in the sense of Cauchy principal value. Let
$$f(t)=y^{(1)}(t).$$
$f(t)$ is a bell-shaped function which is in $L^1$, hence its Fourier transform $F(\omega)$ exists. The existence of the Hilbert transform necessitates that $f(t)$ be in $L^p$ for $1<p<\infty$; for $L^1$ functions, it exists in $L^{1,weak}$ (\cite{Stein1971}, Lemma V.2.8).

We first prove that the odd  derivatives of an odd  sigmoidal function $y(t)$  reach their global extreme values at $t=0$ (Propositions \ref{prop1-boundness}, \ref{prop2-monotone} and  Corollary \ref{corol1-globalextreme}). Then we prove that the global extreme value of the even   derivatives of $y(t)$ is the local extreme values closest to $t=0$ (Proposition \ref{prop3-firstlocakextreme}).
Finally we prove that the sequence of points where the even derivatives of $y(t)$ reach their global extreme values converge to $t=0$ (Corollary \ref{corol2-limitpoint}).
We start by giving a result based solely on the assumption that $F(\omega)$ is positive.

\begin{proposition}\label{prop1-boundness}
Let $f(t)$ be a real, even function such that the Fourier transform   ${\cal F}(f)=F(\omega)$ exists and has constant signature.
Then
 $|f(t)|\le |f(0)|$.
\end{proposition}
\begin{proof}
 We start with the definition and use the inequalities $|\int f |\le \int |f|$ to get
\begin{eqnarray}
|f(t)|&=&\frac{1}{\sqrt{2\pi}}\left|\int_{-\infty}^{\infty} F(\omega)e^{i\omega t}d\omega\ \right|\cr
&\le & \frac{1}{\sqrt{2\pi}}\int_{-\infty}^{\infty} \left|F(\omega)e^{i\omega t}\right|d\omega\ \cr
&= & \frac{1}{\sqrt{2\pi}}\int_{-\infty}^{\infty} \left|F(\omega)\right|d\omega \, .\ \cr
\nonumber
\end{eqnarray}
Since $f(t)$ is symmetric,  $F(\omega)$ is real; furthermore if it has constant signature, then  either
\mbox{$\left|F(\omega)\right|= F(\omega)$} or \mbox{$\left|F(\omega)\right|= -F(\omega)$} for all $\omega$. Thus we have an upper bound for $f(t)$, as

\[
|f(t)| \le
\left\{
\begin{array}{lcl}
  \frac{1}{\sqrt{2\pi}}\int_{-\infty}^{\infty}  F(\omega) \, d\omega= f(0) & , & \left|F(\omega)\right|= F(\omega) \\
 & & \\
  \frac{1}{\sqrt{2\pi}}\int_{-\infty}^{\infty} - F(\omega) \, d\omega=- f(0) &, & \left|F(\omega)\right|=- F(\omega)
\end{array}
\right.
\]
which gives the upper bound $|f(t)|\le |f(0)|.$
\end{proof}

We prove the same result by requiring  the  ``monotonicity" of the envelope of  $y^{(n)}(t)$, expressed in terms of its analytic representation.

\begin{proposition}\label{prop2-monotone}
 Let $f(t)$ be a real, even function such that the magnitude of its analytical representation $|f_A(t)|$ is monotone decreasing on $[0,\infty)$. Then $|f(t)|\le |f(0)|$.
\end{proposition}
\begin{proof}
Recall that $f_A(t)=f(t)+ i  f_h(t)$, where $f_h(t)$ is the Hilbert transform of $f(t)$.
Since $f_A^2(t)=f^2(t)+f_h^2(t)$,  $|f(t)|\le |f_A(t)|$.  When $f(t)$ is even, $f_h(t)$ is odd, hence, $f(0)=f_A(0)$.  Since $|f_A(t)|$ is monotone decreasing, $|f(t)| \le |f_A(0)|=|f(0)|$.
\end{proof}

Note that if the Fourier transform of the first derivative of a sigmoidal curve has constant signature, than all the Fourier transform of all
of its odd derivatives will have constant signature on $[0,\infty)$, hence the local extremum at $t=0$ will be a global one, i.e,
$$\left| y^{(2k+1)}(t) \right|\le \left| y^{(2k+1)}(0) \right|,\quad k\ge 0.$$
We state this below.

\begin{corollary}\label{corol1-globalextreme}
 Let $y(t)$ be an odd sigmoidal function. If $|y_A^{(2k+1)}(t)|$  is monotone decreasing on $[0,\infty)$,
then  $y^{(2k+1)}(t)$ reaches its global extreme value at $t=0$.
\end{corollary}

We will now prove that the global extremum of the even derivatives of an odd sigmoidal curve is the one closest to $t=0$.

 \begin{proposition}\label{prop3-firstlocakextreme}
 Let $f(t)$ be a real, odd function such that the magnitude of its analytical representation $|f_A(t)|$ is monotone decreasing on $[0,\infty)$. Then $|f(t)|\le |f(t_0)|$ where $t_0$ is the location of the first local extremum of $f(t)$ for $t>0$.
\end{proposition}
\begin{proof}
Since $f(t)$ is odd, $f(0)=0$. For simplicity assume that $f'(0)>0$.   Let $t_1$ be the first zero of $f(t)$ for $t>0$. Remark \ref{remark-01-alterne-roots} says that there is a point $t_2<t_1$ at which $f(t_2)=|f_A(t_2)|$ and since  $|f_A(t)|$ is decreasing, $f(t)$ is reaching its first local maximum at some $t_0<t_2$. It follows that $|f(t)|<|f(t_0)|$.
\end{proof}

Finally we will prove that the global extreme values converge to $t=0$.

\begin{corollary}\label{corol2-limitpoint}
Let $y(t)$ be an odd sigmoidal function and assume that
$\left|y_A^{(2k)}(t)\right| $ is monotone decreasing on $[0,\infty)$, then $t=0$ is  the limit point of the global extreme values
of   $ y^{(2k)}(t)$.
\end{corollary}
\begin{proof}
We will prove that the  global extremum  of $y^{(2k+2)}$ occurs earlier than the global extremum   of $y^{(2k)}$.
By Proposition \ref{prop3-firstlocakextreme} above, the global extremum is the first local extremum.  Let $y^{(2k)}(0)=0$, assume that $y^{(2k+1)}(0)>0$ and let $t_1$ be the first intersection of $y^{(2k)}$ with its envelope, $y^{(2k)}(t_1)=\left|y_A^{(2k)}(t_1)\right|$. The global extreme value of $y^{(2k)}$ is at some $t_2<t_1$, since the envelope is decreasing.  Then, $y^{(2k+1)}(t_2)=0$, $y_h^{(2k+1)}(t_2)$ is tangent to the envelope, hence it has its global extremum at $t_3<t_2$.  Finally,  $y_h^{(2k+2)}(t_3)=0$, hence $y^{(2k+2)}(t)$ is tangent to the envelope at this point and it has its global extremum at some $t_4<t_3$.  Hence the global extreme values form a decreasing sequence that converge to $t=0$.
\end{proof}

\subsection{Sigmoidal curves with no symmetry}

Let's now consider a sigmoidal curve $y(t)$ with no symmetry and let $f(t)=y^{(1)}(t)$.
Let $t_m$ be the zero of the first derivative of $f(t)$  and  $t_a$ and $t_b$ be the zeros of its second derivative, with $t_a<t_m<t_b$.
By shifting time, we may assume that $t_m=0$ and
separate the even and odd parts of $f(t)$ as
$$f(t)=f_e(t)+ f_o(t).$$
Then, $$F(\omega)=F_e(\omega)+iF_o(\omega),$$
where $F_e(\omega)$ and $F_o(\omega)$ are real.
We can compute  $f^{(n)}(t)$, $F^{(n)}(\omega)=(i\omega)^n F(\omega) $ and $f_A^{(n)}(t)$ as before. The proof of Proposition \ref{prop3-firstlocakextreme} can be adopted to this case as follows.

 \begin{proposition}\label{prop4-firstlocakextreme}
  Let $f(t)$ be a real function such that the magnitude of its analytical representation $|f_A(t)|$ is monotone decreasing on the half line $[t_*,\infty)$. Then $|f(t)|\le |f(t_0)|$ where $t_0$ is the location of the first local extremum of $f(t)$ for $t>t_*$.
\end{proposition}
\begin{proof}
 For simplicity assume that $f(t_*)\ge 0$. Let $f'(t_*)>0$.   Let $t_1$ be the first zero of $f(t)$ for $t>t_*$ and let $t_2<t_1$ be the point at which $f(t_2)=|f_A(t_2)|$.  Since  $|f_A(t)|$ is decreasing, $f(t)$ is reaching its first local maximum at some $t_0<t_2$. It follows that $|f(t)|<|f(t_0)|$.  If $f'(t_*)<0$ and $f(t_*)>0$ then the global extreme value of $f(t)$ occurs at a point $t_0<t_*$.
\end{proof}

Thus the local extremum of $f(t)$ that is closest to the maximum of $|f_A(t)|$ is its global extremum.  We need to prove two more things: (1) The local extremum should get close to $t_*$ so that there is no gap; (2) The points $t_*$ should converge.

Let $f(t)$ be a bell-shaped curve with Fourier transform $F(\omega)$.  Assume for the time being that $f(t)$ is a small perturbation of an even bell-shaped curve $f_e(t)$ and let $f(t)=f_e(t)+f_o(t)$.  We want to recognize this even-odd decomposition invariantly.  The Fourier transform of $f(t)$ is
$F(\omega)=F_e(\omega)+iF_o(\omega)$ where $F_e(\omega)$ is real and even and $F_o(\omega)$ is real and odd.  But this splitting is not very meaningful, because the shift of a perfectly symmetrical bell-shaped function to $t_o$  will result in a phase of $e^{i t_o \omega}$.  Thus a linear phase corresponds to a shift of time, not to an intrinsic asymmetry. Therefore, if the limit as $\omega\to \pm \infty$ exists, then we can say that $f(t)$ is {\it located where it is supposed to be} and the phase reflects the asymmetry.
We will thus assume that the phase of $F(\omega)$ has a limit as $\omega\to \infty$.  Then, the inverse Fourier transform of the real part of $F(\omega)$
will be the even part of $f(t)$ and the even and odd parts are invariantly defined.  Note that the even part reaches its maximum at $t=0$, but the maximum of $f(t)$ is at some $t_m$ that lies at the right or at the left of $t=0$ depending on the sign of the odd part.

We have seen that the Fourier transform of the higher derivatives $F^{(n)}(\omega)=i^n \omega^n F(\omega) $ could be approximated by $i^n F(\omega-\omega_n)+ i^n F(\omega+\omega_n)$ where $\pm\omega_n$ are the points where $|F^{(n)}(\omega)|$ reaches its maximum.  Note that up to the factor of $i^n$ the phase of $F^{(n)}(\omega)$ and $F(\omega)$ are the same. Again for the time being assume that $n$ is even  consider $\pm f^{(n)}(t)$. Note that $F(0)$ is real whenever $f(t)$ is real. Thus, the time shift in $\pm f^{(n)}(t)$ is the phase of $F(\omega)$ at $\omega_n$, divided by $\omega_n$.  If the phase of $F(\omega)$ has a limit as $\omega\to \infty$, then the time shift always tends to zero. Thus, $t=0$ is the critical point for the asymmetrical case also, provided that the initial time shift is arranged so that the even part of $f(t)$ is centered at $t=0$.
For the odd derivatives, the same scheme applies to the Hilbert transform $f_h(t)$.  This is the reason why the even derivatives converge to zero faster than the odd derivatives.

\vskip 0.2cm
\noindent
{\bf Assumption:}
  $F^{(n)}(\omega)$ has only two maxima at $\pm\omega_n$ and    $|f_A^{(n)}(t)|$ has a single local maximum  at $t_n$.
   \vskip 0.2cm

   As before, we can prove that the global maximum of $|y^{(n+1)}(t)|$ is at $t_n$ if $n+1$ is odd, and it is the local extremeum closest to $t_n$ if $n+1$ is even.

The localization condition is now expressed as $|(i\omega)^n F(\omega)|$ having a single maximum at $\omega_n$ for $\omega>0 $.  Now, as opposed to the case, as $F(\omega)$ is not real, $(i\omega)^n F(\omega)$ has a phase $\phi(\omega)$, i.e,
$$(i\omega)^n F(\omega)= |\omega^n F(\omega)| e^{i\phi(\omega)}.$$
We add and subtract a linear phase factor $\alpha_n\omega$ to get
$$(i\omega)^n F(\omega)= |\omega^n F(\omega)| e^{i\phi(\omega)+i\alpha_n\omega} \ e^{-i\alpha_n \omega}.$$
If we choose
$$\alpha_n=-\phi(\omega_n)/\omega_n$$
then the phase at $\omega_n$ will be exactly $-\alpha_n\omega_n$.  Finally using the time shift property, we get
$${\cal F}^{-1} \left( (i\omega)^n F(\omega)\right) = y^{(n+1)}(t-\alpha_n).$$

\section{The Standard and Generalized Logistic Growth}

The standard logistic growth curve is the solution of the differential equation
\[y'=1-y^2,\quad y(0)=0.\]

This equation can be solved as $y(t)=\tanh(t)$ and its first derivative, $y'(t)={\rm sech}^2(t)$
is the well known $1$-soliton solution of the Korteweg-deVries equation.

\begin{figure}[h]
    {\centering
    $\begin{array}{cc}
       \includegraphics[width=0.5\textwidth]{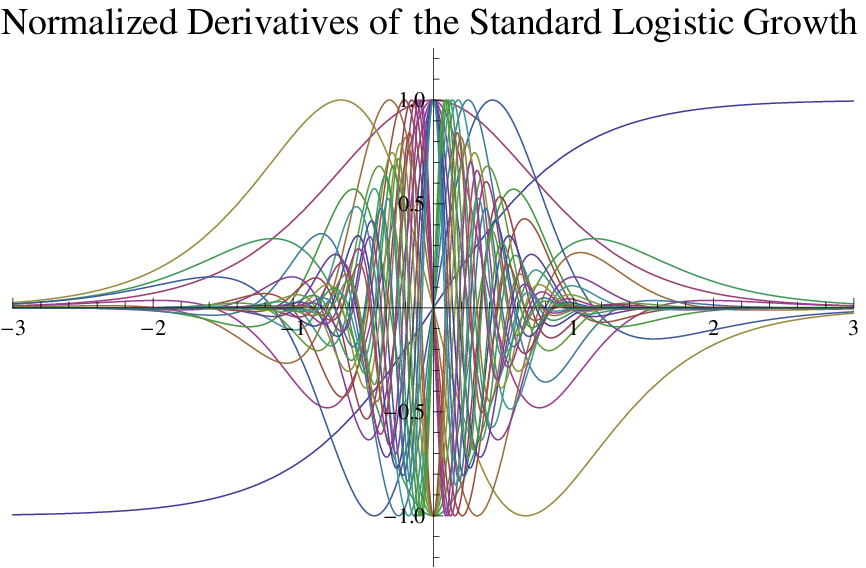}
       & \includegraphics[width=0.5\textwidth]{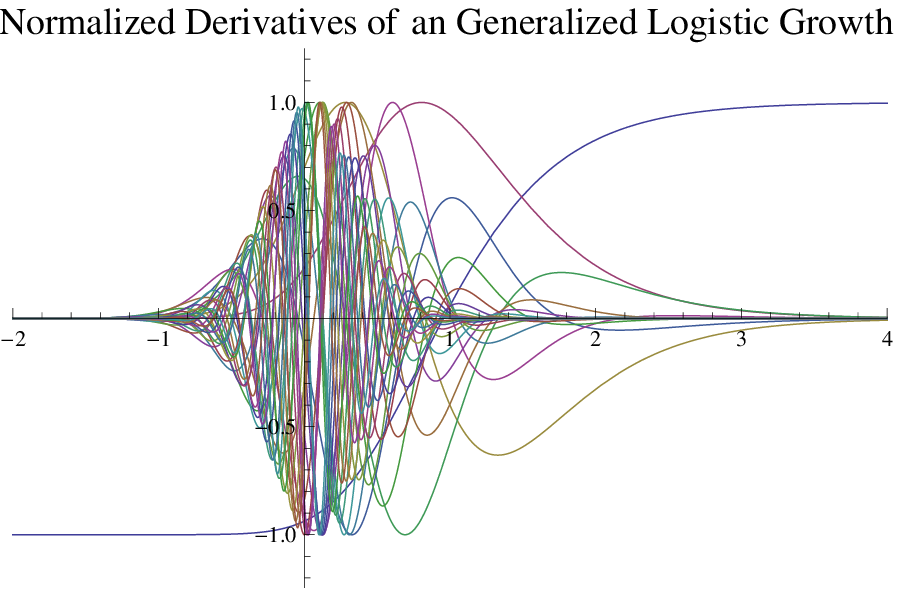}
       \\
       \rm (a) & \rm (b)
       \\
    \end{array}$}
          \caption{Normalized derivatives of the standard logistic growth (a) and normalized derivatives of the generalized logistic growth ($\beta=1, k=1, \nu=1/5$) (b) up to order $30$.\label{figure-30derivatives}}
  \end{figure}

\begin{figure}[h]
    {\centering
    $\begin{array}{cc}
       \includegraphics[width=0.5\textwidth]{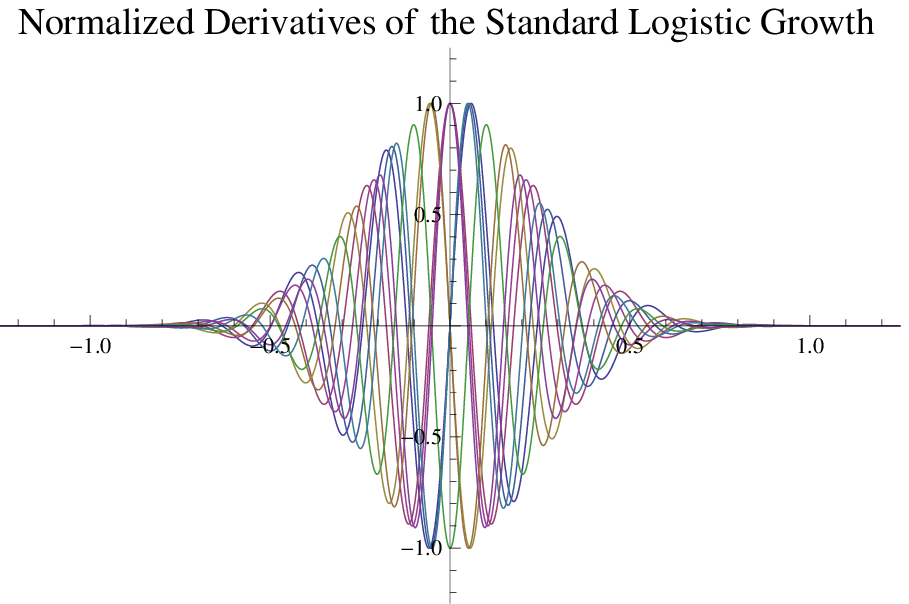}
       & \includegraphics[width=0.5\textwidth]{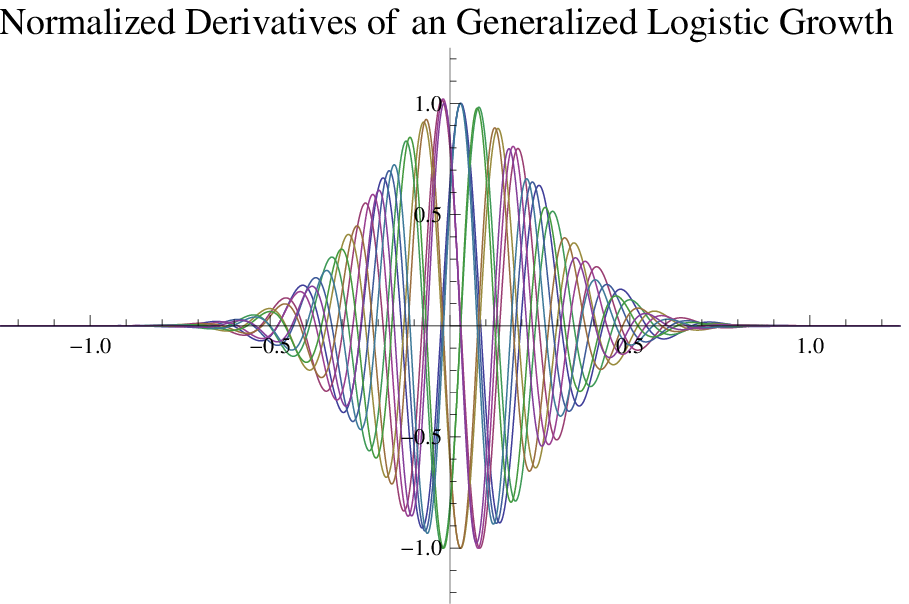}
       \\
       \rm (a) & \rm (b)
       \\
    \end{array}$}
          \caption{Normalized derivatives of the standard logistic growth (a) and normalized derivatives of the generalized logistic growth ($\beta=1, k=1, \nu=1/5$) (b) of order between $40$ and $50$. \label{figure-4050derivatives}}
  \end{figure}

\noindent The generalized logistic growth curve with horizontal asymptotes at $-1$ and $1$ is given by
\begin{equation} \label{asym-log-curve}
y(t)=-1+\frac{2}{\left[1+ke^{-\beta t}\right]^{1/\nu}},
\end{equation}
where   $k>0$, $\beta>0$ and $\nu>0$. The parameter $k$ can be adjusted by a time shift, $\beta$ corresponds to a scaling of time; $\nu$ is the key parameter that determines the shape of the growth, the asymmetry depending  on whether it is less than or greater than $1$.
This formula reduces to the standard logistic growth for $\nu=1$, $k=1$, $\beta=2$. We are plotting the normalized derivatives of the standard and  generalized logistic growth functions of order up to $30$ in Figure~\ref{figure-30derivatives} and of order between $40$ and $50$ in Figure~\ref{figure-4050derivatives}.

\noindent \textbf{Fourier transform of the logistic growth:}
 Let $f(t)={\rm sech}^2(t)$.  Using the integral formula
$$\int_{-\infty}^{\infty} \frac{e^{ixy}}{\cosh^2(x)} dx= \frac{\pi y}{\sinh\left(\frac{\pi y}{2}\right)},$$
the Fourier transform of $f(t)$ is found as below:
$$ f(t)={\rm sech}^2(t)\quad \to \quad F(\omega)=\sqrt{\frac{2}{\pi}}\frac{(\pi\omega/2)}{\sinh(\pi\omega/2) }.$$
Thus, the Fourier transform pairs for the derivatives are
$$ f^{(n)}(t)=\frac{d^n}{dt^n}\large({\rm sech}^2(t)\large)\quad \to \quad
F^{(n)}(\omega)=(i)^n\sqrt{\frac{2}{\pi}}    \left({\frac{2}{\pi}}\right)^n    \frac{(\pi\omega/2)^{n+1}}{\sinh(\pi\omega/2) }.$$

\begin{figure}[h]
    {\centering
    $\begin{array}{cc}
       \includegraphics[height=5cm]{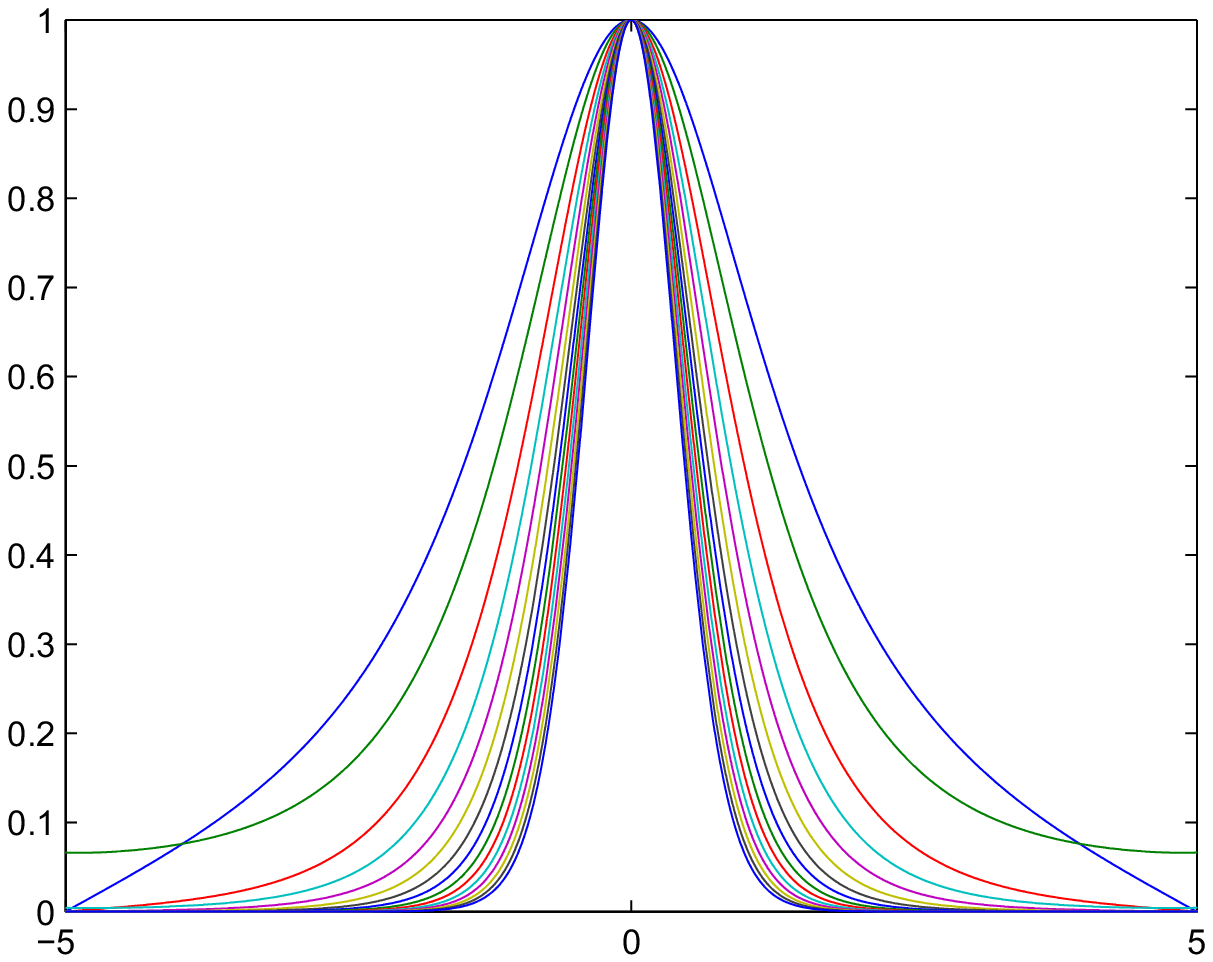}
       & \includegraphics[height=5cm]{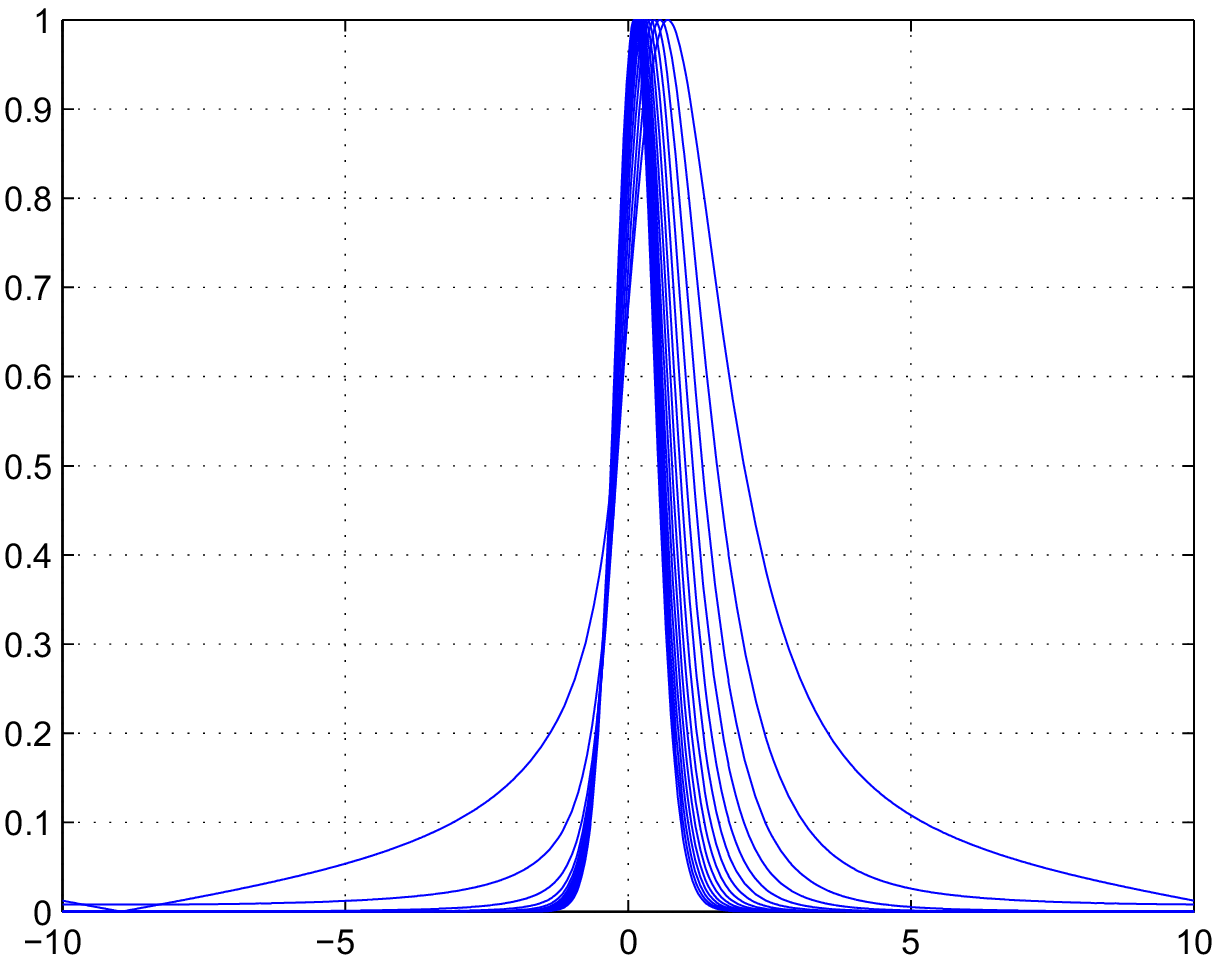}
       \\
       \rm (a) & \rm (b)
       \\
    \end{array}$}
     \caption{ Normalized envelopes of the derivatives of the standard logistic growth function  up to order $15$ (a), Normalized envelopes of the derivatives of the generalized logistic growth function ($\beta=2$, $k=1$, $\nu=1/5$) (b) up to order $15$.\label{figure-envelopeofderivatives}}
\end{figure}

Since the Fourier transform is positive the assumption of Proposition 1 holds.  Unfortunately it is not possible to prove that Hilbert transforms are monotone. We present normalized envelopes in Figure~\ref{figure-envelopeofderivatives} where wee see the monotonicity of the envelopes. Also one can observe from the figure the symmetry and asymmetry properties in the standard and generalized cases respectively.

\begin{figure}[h]
    {\centering
    $\begin{array}{cc}
       \includegraphics[width=0.5\textwidth]{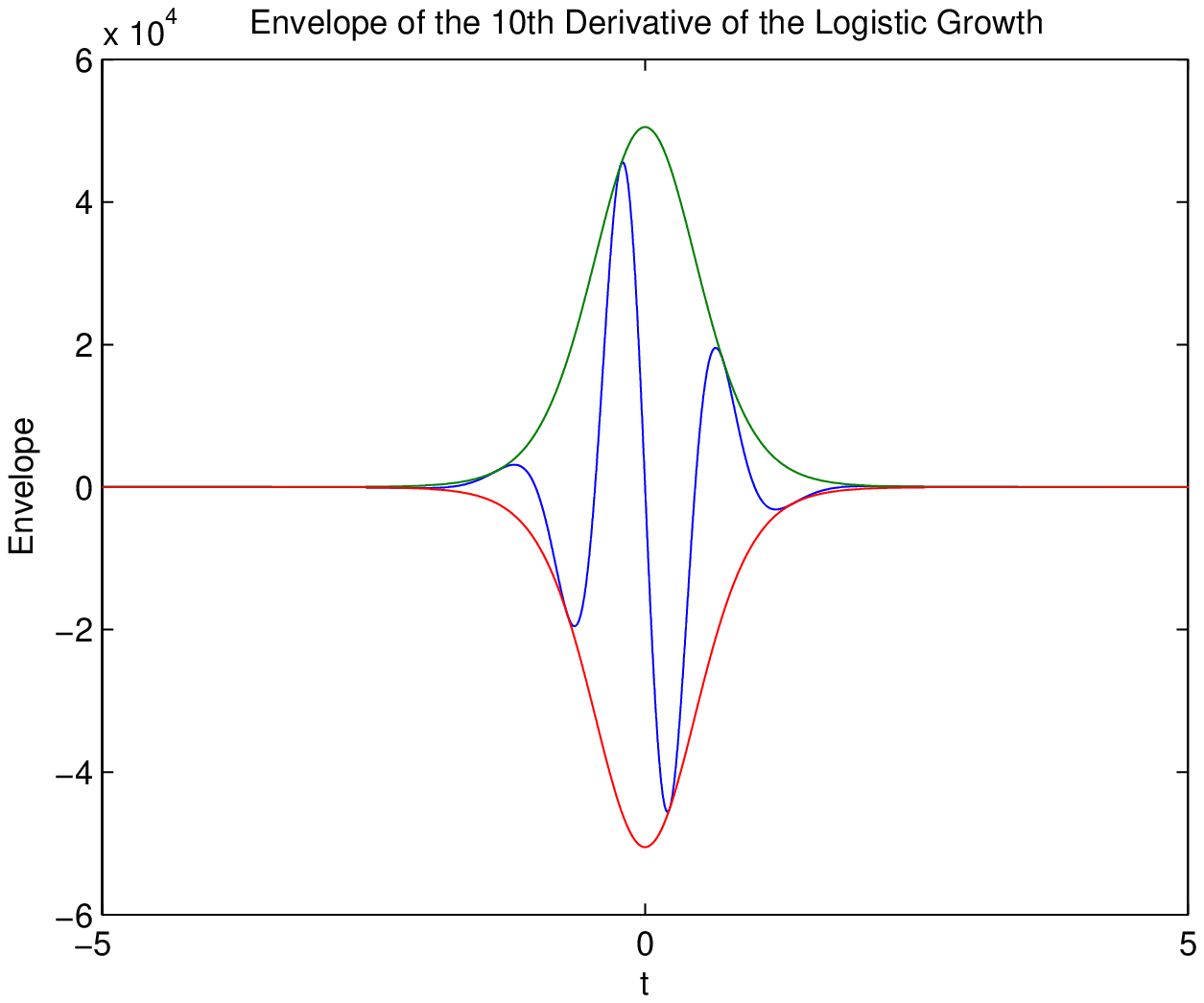}
       & \includegraphics[width=0.5\textwidth]{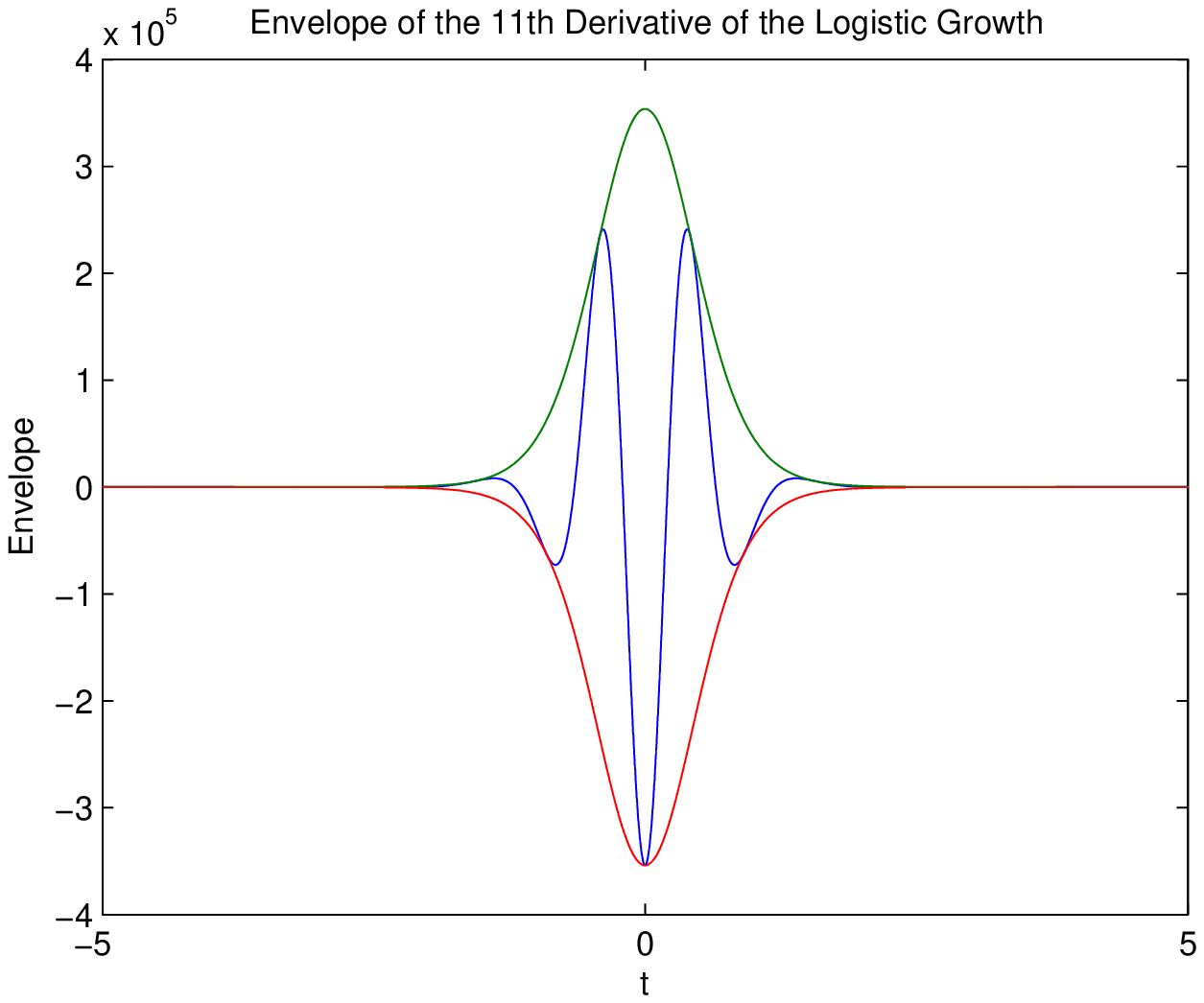}
       \\
       \rm (a) & \rm (b)
       \\
    \end{array}$}
    \caption{ The envelope of the $10$th derivative of the standard logistic growth function (a), The envelope of the $11$th derivative of the standard logistic growth function (b).\label{figure-envelope-10-11}}
\end{figure}

\begin{figure}[h]
    {\centering
    \includegraphics[width=0.6\textwidth]{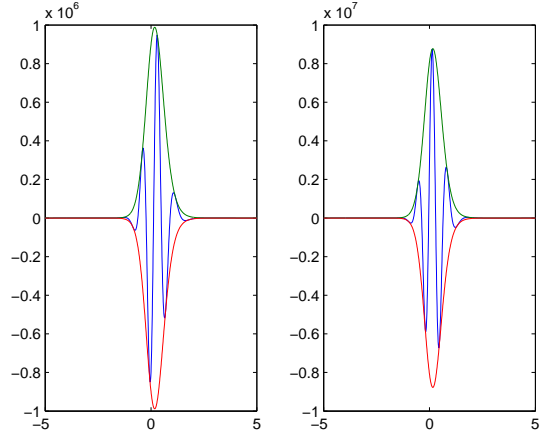}
    \caption{ The envelope of the $10$th derivative of the generalized logistic growth function (a), The envelope of the $11$th derivative of the generalized  logistic growth function (b).\label{figure-asymmetricenvelope-10-11}}}
\end{figure}

In Figure~\ref{figure-envelope-10-11} (and Figure~\ref{figure-asymmetricenvelope-10-11}) we give more detailed figures include the envelopes of the  $10$th and $11$th derivatives of the standard logistic growth function (and the generalized logistic growth function).


\begin{figure}[h]
{\centering
\includegraphics[width=0.75\textwidth]{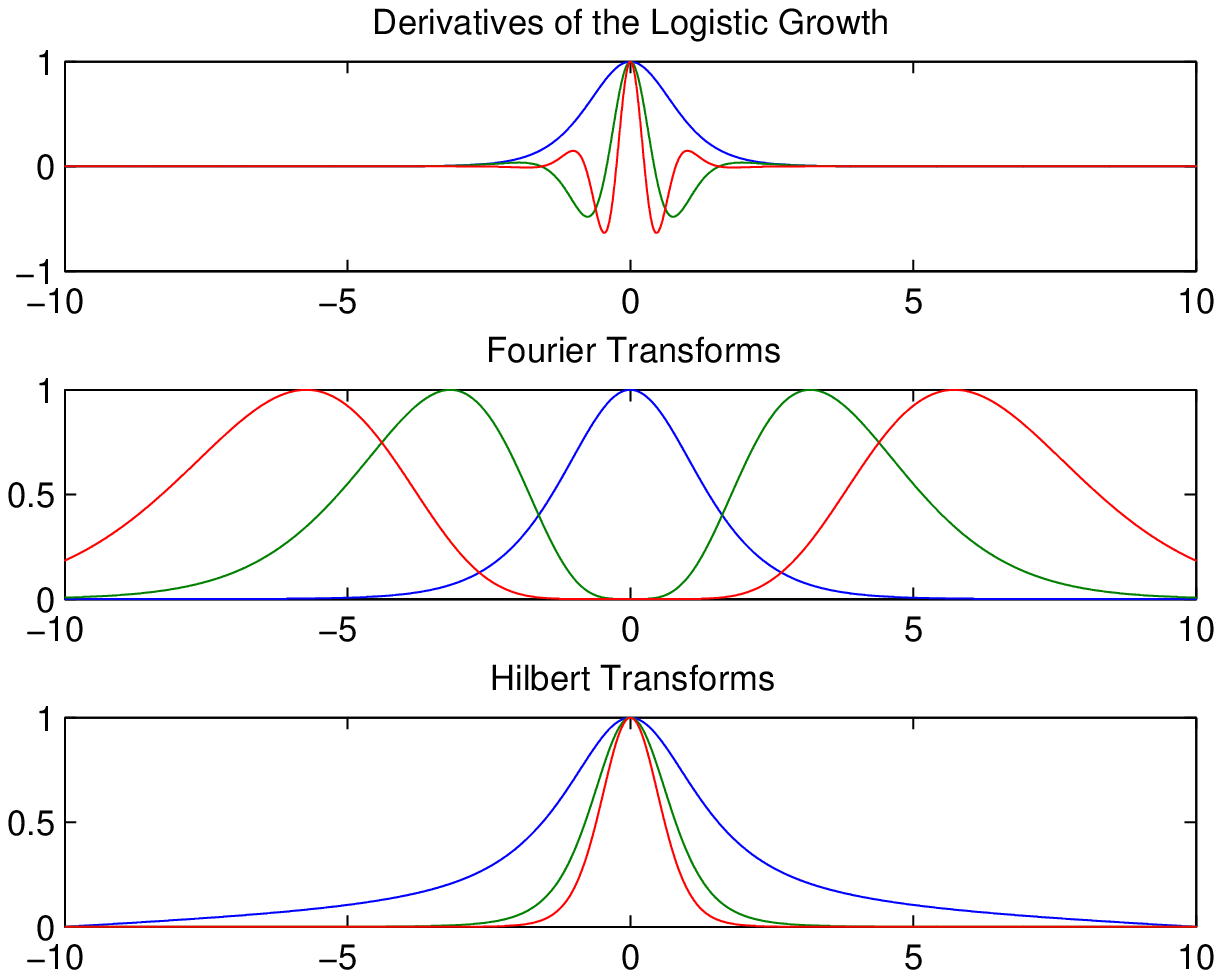}
\caption{Normalized  graphs of the time domain variations, the magnitude of the Fourier transform and the magnitude of the analytic representation   for $y(t)$, $y^{(4)}(t)$ and  $y^{(8)}(t)$. \label{figure-y(t)-4-8} }}
\end{figure}

In Figure~\ref{figure-y(t)-4-8}, we present respectively the graphs of $y(t)$, $y^{(4)}(t)$ and $y^{(8)}(t)$, the magnitude of their Fourier transform and the magnitude of their analytic representation.

\begin{figure}
    {\centering
    \includegraphics[width=0.7\textwidth]{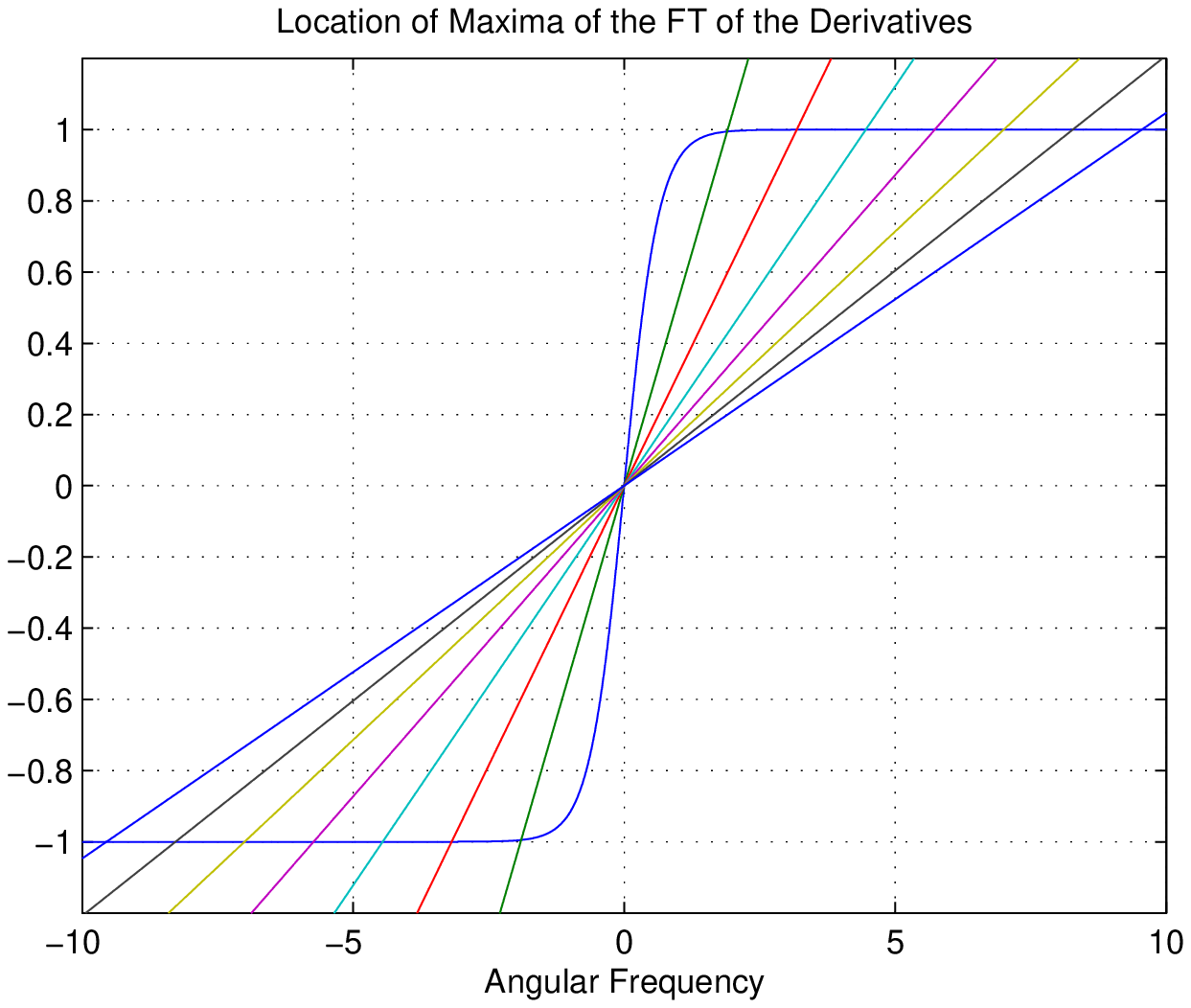}
    \caption{ The Fourier transform of each derivative has a single maximum for positive frequencies. \label{figuremaxima}}}
\end{figure}

From  Figure~\ref{figure-y(t)-4-8} we see that the Fourier transform of higher order derivatives are localized humps moving to $\pm \infty$. In fact we can prove that the Fourier transform of the $n$th derivative has a single maximum for positive frequencies.  This can be seen by taking the derivative of the transform and setting equal to zero (see Figure \ref{figuremaxima}).

\begin{figure}
    {\centering
    \includegraphics[width=0.75\textwidth]{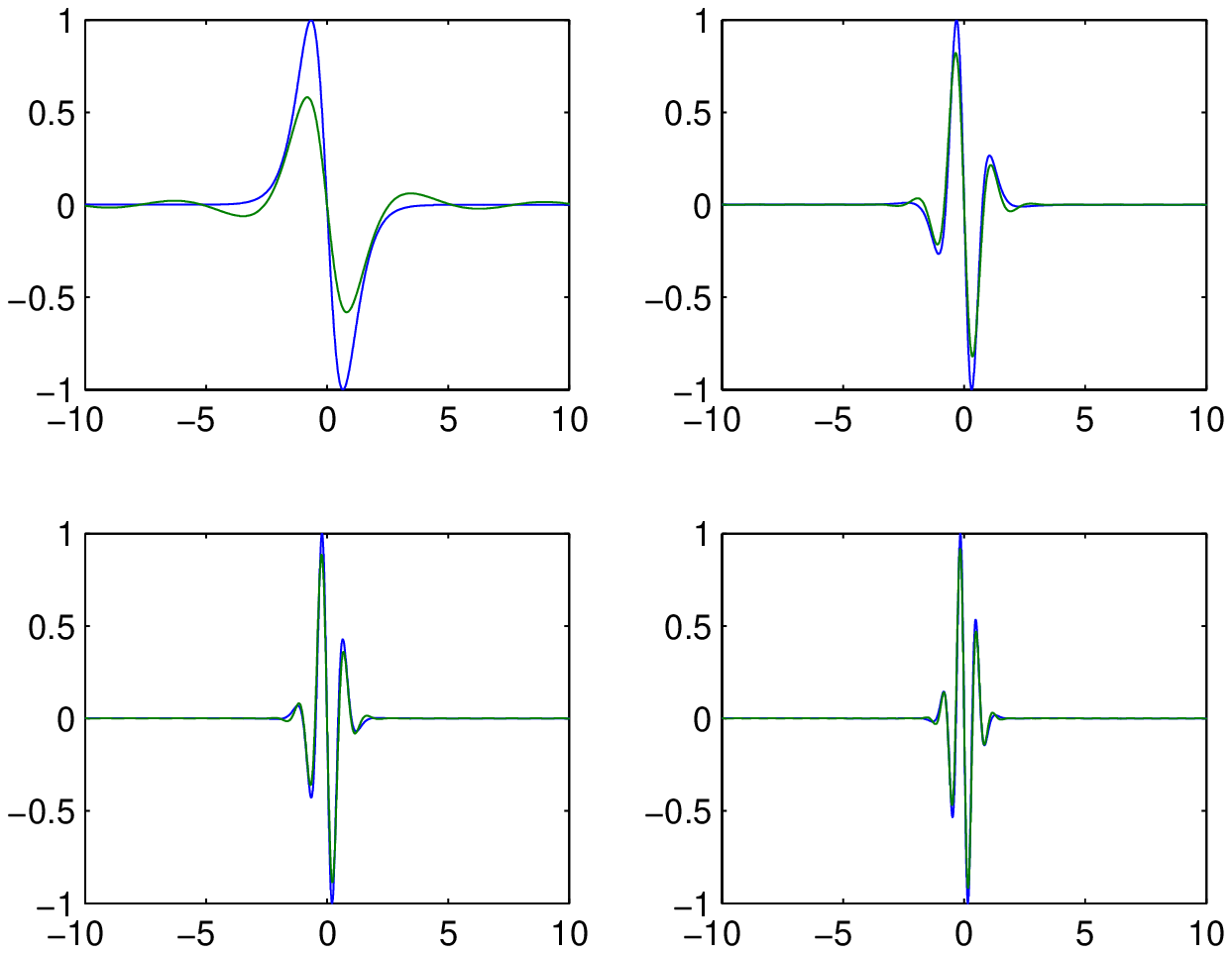}
    \caption{In the time domain, the derivatives look like a modulated waveform. \label{figuretimedomain}}}
\end{figure}

In the time domain, the $n$th derivative looks in fact like a modulated waveform (see Figure~\ref{figuretimedomain}).

\textbf{Discussion of the Assumptions:}
In the standard case, for  $n$ even,  we can use Proposition \ref{prop1-boundness}  to conclude that $t=0$ is the global extremum provided that $F^{(n)}(\omega)$ is positive.  Since this is satisfied for $F(\omega)$ above, we conclude that all odd derivatives of the standard logistic growth function have their global extremum at $t=0$. For odd derivatives of $f(t)$, hence even derivatives of the standard logistic growth (notice that $f(t)=y^{(1)}(t)$), we  need to compute the Hilbert transform of $f^{(2k+1)}(t)$. Unfortunately this is not possible to do analytically.

For the standard logistic growth function we can easily see that all  assumptions are satisfied.

As indicated before, the Fourier transform of the generalized logistic growth function is presented in Appendix. In order to compute the phase of $F(\omega)$ for this case we assume that $1/\nu$ is an integer, take $1/\nu=n$, then we obtain
\[
F(\omega)=\sqrt{\frac{2}{\pi}}\, e^{-\frac{i\omega}{\beta} \ln k} \left[(n-1)-\frac{i\omega}{\beta} \right] \left[(n-2)-\frac{i\omega}{\beta} \right]\cdots \left[1-\frac{i\omega}{\beta} \right]\left[-\frac{i\omega}{\beta} \right]\left[\frac{i\omega}{\beta} \right].
\]
by using Eq.(\ref{asym-log-curve-gamma-Fw}) and applying the property given by Eq.(\ref{asym-log-curve-gamma}) repeatedly (see Appendix).
For simplicity we take $k=1$. Using the fact that $\Gamma(\overline{z})=\overline{\Gamma(z)}$, we can see that the phase of $F(\omega)$ is
\[
\phi_n(\omega) \tan^{-1} \left( \frac{-\omega/\beta}{n-1} \right) \cdots \tan^{-1} \left( \frac{-\omega/\beta}{2} \right)\tan^{-1} \left( \frac{-\omega/\beta}{1} \right).
\]
It then follows that the phase factor $\alpha_n = {\rm tan^{-1}} (\omega_n/\beta)/ \omega_n$, hence $\alpha_n$ converges to zero.  This agrees with the observations.

\section*{Acknowledgments}
This work is partially supported by The Scientific and Technological Research Council of Turkey (TUBITAK).

\section*{Appendix}

As given in Section 4, Eq.(\ref{asym-log-curve}), the generalized logistic growth curve with horizontal asymptotes at $-1$ and $1$ is given by
\begin{equation*}
y(t)=-1+\frac{2}{\left[1+ke^{-\beta t}\right]^{1/\nu}},
\end{equation*}
where   $k>0$, $\beta>0$ and $\nu>0$. Then the first derivative of the generalized logistic growth curve is
\begin{equation*} \label{asym-log-curve-derivation}
f(t)=y^{(1)}(t)=\frac{2k\beta}{\nu}\left[1+ke^{-\beta t}\right]^{-\frac{1}{\nu}-1} e^{-\beta t}.
\end{equation*}
and its Fourier transform is given by
\begin{equation*} \label{asym-log-curve-fourier}
F(\omega)=\frac{1}{\sqrt{2\pi}} \frac{2k\beta}{\nu}
\int_{-\infty}^{\infty}  e^{-i\omega t}
e^{-\beta t }\left[1+ke^{-\beta t}\right]^{-\frac{1}{\nu}-1}\ dt.
\end{equation*}

The indefinite integral is evaluated in terms of the hypergeometric functions as below. Let
$$I(\omega)= \int_{-\infty}^{\infty}  e^{-i\omega t}
e^{-\beta t }\left[1+ke^{-\beta t}\right]^{-\frac{1}{\nu}-1}\ dt.  $$
Then
$$I(\omega)=
\frac{i\nu}{k(\nu \omega+i\beta )} e^{-i\omega t}
\left(1+\frac{e^{bt}}{k}\right)^{\frac{1}{\nu}}(ke^{-\beta t}+1)^{-\frac{1}{\nu}}
\ _2F_1\left(1+\frac{1}{\nu},\frac{1}{\nu}-\frac{i\omega}{\beta},
\frac{1}{\nu}-\frac{i\omega}{\beta}+1,-\frac{e^{\beta t}}{k}\right).
  $$
The integral can be simplified as
$$I(\omega)=
\frac{i\nu}{k(\nu \omega+ib)}
e^{-i\omega t}
\left(\frac{e^{\beta t}}{k}\right)^{\frac{1}{\nu}}
\ _2F_1\left(
1+\frac{1}{\nu},
\frac{1}{\nu}-\frac{i\omega}{\beta},
\frac{1}{\nu}-\frac{i\omega}{\beta}+1,
-\frac{e^{\beta t}}{k}\right).
  $$
The hypergeometric function   has the property
$$\ _2F_1(a,b,c,0)=0.$$
Since $e^{-i\omega t}$ remains bounded and $\beta >0$, the indefinite integral vanishes as $t\to -\infty$.
It follows that the Fourier transform is given by
$$F(\omega)=\frac{1}{\sqrt{2\pi}}
\frac{2k\beta}{\nu}
\frac{i\nu}{k(\nu \omega+i\beta )}
\lim_{t\to \infty }
\left[
e^{-i\omega t}
\left(\frac{e^{\beta t}}{k}\right)^{\frac{1}{\nu}}
\ _2F_1\left(
1+\frac{1}{\nu},
\frac{1}{\nu}-\frac{i\omega}{\beta},
\frac{1}{\nu}-\frac{i\omega}{\beta}+1,
-\frac{e^{\beta t}}{k}\right)\right].$$
After simplifications we get
$$F(\omega)=\sqrt{ \frac{2}{\pi}  }
\frac{i\beta}{\nu \omega+i\beta}
\lim_{t\to \infty }
\left[
e^{-i\omega t}
\left(\frac{e^{\beta t}}{k}\right)^{\frac{1}{\nu}}
\ _2F_1\left(
1+\frac{1}{\nu},
\frac{1}{\nu}-\frac{i\omega}{\beta},
\frac{1}{\nu}-\frac{i\omega}{\beta}+1,
-\frac{e^{\beta t}}{k}\right)\right].$$
Defining
$$z=-\frac{e^{\beta t}}{k},$$
and using
$$e^{-i\omega t}= \left( e^{\beta t }\right)^{-\frac {i\omega}{\beta}}
=\left( -k z\right)^{-\frac {i\omega}{\beta}}
= \left(k \right)^{-\frac {i\omega}{\beta}}
\left( - z\right)^{-\frac {i\omega}{\beta}},$$
we get
$$F(\omega)=\sqrt{ \frac{2}{\pi}  }
\frac {i\beta  \left(k \right)^{-\frac {i\omega}{\beta}}     }
      {\nu \omega+i\beta                                     }
\lim_{t\to \infty }
\left[
\left(-z\right)^{-\frac {i\omega}{\beta} +\frac{1}{\nu}       }
\ _2F_1\left(
1+\frac{1}{\nu},
\frac{1}{\nu}-\frac{i\omega}{\beta},
\frac{1}{\nu}-\frac{i\omega}{\beta}+1,
z\right)\right].$$
This expression is  of the form
$$F(\omega)=\sqrt{ \frac{2}{\pi}  }
\frac {i\beta  \left(k \right)^{-\frac {i\omega}{\beta}}     }
      {\nu \omega+i\beta                                     }
 \times \lim_{z\to -\infty }
\left[(-z)^b \ _2F_1(a,b,b+1,z)\right],$$
where
$$a=1+\frac{1}{\nu},\quad b=\frac{1}{\nu}-\frac{i\omega}{\beta}.$$
We use the linear transformation properties of the hypergeometric functions given in \cite{Abrom}
\begin{eqnarray}
\ _2F_1(a,b,c,z)&=&
\frac{\Gamma(c)\Gamma(b-a)}{\Gamma(b)\Gamma(c-a)}(-z)^{-a}
  \ _2F_1(a,1-c+a,1-b+a,1/z)\nonumber\cr
&&+\frac{\Gamma(c)\Gamma(a-b)}{\Gamma(a)\Gamma(c-b)}(-z)^{-b}
  \ _2F_1(b,1-c+b,1-a+b,1/z).\nonumber
  \end{eqnarray}
For $c=b+1$, we get
\begin{eqnarray}
\ _2F_1(a,b,b+1,z)&=&
\frac{\Gamma(b+1)\Gamma(b-a)}{\Gamma(b)\Gamma(1+b-a)}(-z)^{-a}
  \ _2F_1(a,-b+a,1-b+a,1/z)\nonumber\cr
&&+\frac{\Gamma(b+1)\Gamma(a-b)}{\Gamma(a)\Gamma(1)}(-z)^{-b}
  \ _2F_1(b,0,1-a+b,1/z).\end{eqnarray}
It is known that whenever $a=0$ or $b=0$,
$\ _2F_1(a,b,c,z)$ reduces to a $1$ (this is best seen in the definition of the hypergeometric function given in  terms of Pochhammer symbols, defined by the product of $n$  terms $(a)_n=a (a+1)(a+2)\cdots (a+n-1)$.
Substituting this in the Fourier transform, we obtain
\begin{eqnarray*}
F(\omega)&=&\sqrt{ \frac{2}{\pi}  }
\frac {i\beta  \left(k \right)^{-\frac {i\omega}{\beta}}     }
      {\nu \omega+i\beta                                     }
      \times
\\
&& \lim_{z\to -\infty }
\left[
\frac{\Gamma(b+1)\Gamma(b-a)}{\Gamma(b)\Gamma(1+b-a)}(-z)^{-a+b}
  \ _2F_1(a,-b+a,1-b+a,1/z)  \right.
\\
&&  \left.
+\frac{\Gamma(b+1)\Gamma(a-b)}{\Gamma(a)\Gamma(1)}
  \ _2F_1(b,0,1-a+b,1/z)\right].
  \end{eqnarray*}
Since the real part of $-a+b$ is negative, and
$\ _2F_1(a,b,c,0)=0$, the first term vanishes in the limit and we get
$$F(\omega)=\sqrt{ \frac{2}{\pi}  }
\frac {i\beta  \left(k \right)^{-\frac {i\omega}{\beta}}     }
      {\nu \omega+i\beta                                     }
\lim_{z\to -\infty }
\left[
\frac{\Gamma(b+1)\Gamma(a-b)}{\Gamma(a)\Gamma(1)}
  \ _2F_1(b,0,1-a+b,1/z)\right],$$
  which finally gives
  $$F(\omega)=\sqrt{ \frac{2}{\pi}  }
\frac {i\beta  \left(k \right)^{-\frac {i\omega}{\beta}}     }
      {\nu \omega+i\beta                                     }
\frac{\Gamma(b+1)\Gamma(a-b)}{\Gamma(a)\Gamma(1)}.$$
Substituting the expressions of $a$ and $b$, using the property
\begin{equation} \label{asym-log-curve-gamma}
\Gamma(1+x)=x\Gamma(x)
\end{equation}
and simplifying the resulting expressions, we obtain
\begin{equation} \label{asym-log-curve-gamma-Fw}
F(\omega)=\sqrt{ \frac{2}{\pi}  }
\frac{ k^{-\frac {i\omega}{\beta}}    }
      {\Gamma(1/\nu)  }
\Gamma\left(\frac{1}{\nu}-\frac{i\omega}{\beta}\right)
\Gamma\left(1+\frac{i\omega}{\beta}\right).
\end{equation}

For $\nu=1$, we use the property
\begin{equation*} \label{asym-log-curve-gamma-sin}
\Gamma(x)\Gamma(1-x)=\frac{\pi}{\sin(\pi x)}
\end{equation*}
together with $\sin(ix)=i\sinh(x)$ to get the Fourier transform of the standard logistic growth function.

\newpage

\noindent{\bf  Ayse Humeyra Bilge}\\
        {\it  Faculty of Engineering and Natural Sciences, Kadir Has University}\\
        {\it  Istanbul, Turkey}\\
        {\it  e-mail: ayse.bilge@khas.edu.tr}\\ {} \\
        {\bf  Yunus Ozdemir}\\
        {\it  Department of Mathematics, Faculty of Science, Anadolu University}\\
        {\it  Eskisehir, Turkey}\\
        {\it  e-mail: yunuso@anadolu.edu.tr}


\begin{thebibliography}{100}

\bibitem{Abrom} M. Abramowitz, I.A. Stegun, Handbook of Mathematical Functions, Dover, New York, 1972.


\bibitem{BP2013} A.H. Bilge and O. Pekcan, {\it A mathematical Description of the Critical Point in Phase Transitions}, Int. J. Mod. Phys. C, vol. 24, 2013.

\bibitem{BP2012} A.H. Bilge, O. Pekcan and V. Gurol, {\it Application of epidemic models to phase transitions},
Phase Transitions, vol. 85, pp. 1009-1017, 2012.

\bibitem{Duoan} J. Duoandikoetxea, Fourier analysis, Graduate Studies in Mathematics, vol. 29, AMS, Providence, RI, 2001.


\bibitem{Grad} I.S. Gradshteyn and I.M. Ryzhik (Eds: A. Jeffrey and D. Zwillinger),  Table of Integrals, Series, and Products,  Elsevier Inc., 2007.


\bibitem{PapoulisFourier} A. Papoulis, The Fourier Integral and its Applications, McGraw-Hill Co., New York, 1962.

\bibitem{Polya1942} G. Polya, {\it On the zeros of the derivatives of a function and its analytic character},
 Bull. Amer. Math. Soc., vol.49, pp. 178-191, 1942.

 \bibitem{IMA1} R.C. Sharpley and V. Vatchev, {\it Analysis of the Intrinsic Mode Functions}, Constr. Approx., vol. 24, pp.17–47, 2006.

 \bibitem{Stein1971} E.M. Stein and G.L. Weiss, Introduction to Fourier Analysis on Euclidean Spaces, Princeton University Press, 1971.


\end{thebibliography}
\end{document}